\documentclass[12pt]{article}
\usepackage{amsmath}
\usepackage{amssymb}
\usepackage{amsthm}
\usepackage{amsfonts}
\usepackage{ulem}
\usepackage{cancel}
\usepackage{graphicx}
\usepackage{dsfont} % for unit operator  (indicator function) \mathds{1}

\usepackage[mathscr]{eucal}

\usepackage{txfonts}

\usepackage{bbm,epsfig,graphics,epic,color,rotating,color}
\numberwithin{equation}{section}

\newcommand{\R}{{\mathbb R}}
\newcommand{\Z}{{\mathbb Z}}

\newcommand{\N}{{\mathfrak{N}}}
\newcommand{\Ns}{{\mathscr{N}}}

\newcommand{\wh}[1]{\widehat{#1}}

\newcommand{\ed}{\mathrm{d}}

\newcommand{\reff}[1]{(\ref{#1})}

\renewcommand{\k}{\varkappa}

\renewcommand\phi{\varphi}

\newcommand{\be}{\begin{equation}}
\newcommand{\ee}{\end{equation}}
\newcommand{\bel}[1]{\begin{equation}\label{#1}}
\newcommand{\bea}{\begin{eqnarray}}
\newcommand{\eea}{\end{eqnarray}}
\newcommand{\balign}{\begin{align}}
\newcommand{\ealign}{\end{align}}
\newcommand{\ba}{\begin{array}}
\newcommand{\ea}{\end{array}}
\newcommand{\bfig}{\begin{figure}}
\newcommand{\efig}{\end{figure}}

\newcommand{\C}{{\mathbb C}}

\newtheorem{theorem}{Theorem}[section]

\theoremstyle{definition}

\title{Large Emission Regime in Mean Field Luminescence}

\author{ E. Pechersky$^{1}$,\and S. Pirogov$^{1}$,\and G. M. Sch\"utz$^{2}$, \and A.~Vladimirov$^{1}$  and  A. Yambartsev$^3$ }
\date{}

\begin{document}

\maketitle

$^1$ Institute for Information Transmission Problems, 19, Bolshoj Karetny, Moscow, 127994, RF

$^2$Interdisziplin\"ares Zentrum f\"ur Komplexe Systeme, Universit\"at Bonn, Br\"uhler Str. 7, 53119 Bonn, Germany

$^3$Institute of Mathematics and Statistics, University of S\~ao Paulo (USP), S\~ao Paulo 05508-090,
SP, Brazil

\begin{abstract}
We study a class of random processes on $N$ particles which can be interpreted as
stochastic model of luminescence. Each particle can stay in one of two states: Excited state 
or ground state. Any particle at ground state is excited with a constant rate (pumping). 
The number of excited particles decreases by means of photon emission through 
interactions of the particles. We analyse the rare event of flashes, i.e., the emission 
of a very large number of photons $B$ during a fixed time interval $T$.
We employ the theory of large deviations to provide the asymptotics of
the probability of such event when the total number of particles
$N$ tends to infinity. This theory gives us also the optimal
trajectory of scaled process corresponding to this event.
The stationary regime of this process we call the large 
emission regime. In several cases we prove that in the large
emission regime a share of excited particles in a system is
stable under the changes of the pumping and emission rates.
\end{abstract}



\section{Introduction}\label{intr} 

We study large fluctuations in a family $\Theta$ of continuous-time Markov processes
which are mean field models of luminescence where $N$ atoms (or, more generally, some
species of particles)
can be either in the ground state
0 or in an excited state 1 and where the transition from the excited state to the ground state
is accompanied by the emission of light. Our interest is in the luminescence distribution, i.e., 
in the probability of the event that the cumulative number of the transitions 
from 1 to 0 during a finite time-interval $T$ is at least $B$.
We use the rate function found by large deviation theory on the path level  to evaluate 
the asymptotics of the probability of some functionals that quantify the asymptotic properties
of the luminescence distribution when $B$, $N$
are large.
In addition we are interested in the asymptotic  
behaviour of some functionals of the process trajectories under the large fluctuations.

We consider several scenarios of luminescence. In our setting all possible configurations form 
the state space $\mathbb{C}=\{0,1\}^N$, and the Markov processes $\Theta$ that we consider 
describe the random transitions within the state space $\mathbb{C}$, when a single particle or a 
group of particles changes its state. 
In order to change their states the particles interact with each other.  
To this end the particles create 
groups of some specific size with some specific relations of their states within the group
as defined below.
For any Markov process $\theta\in\Theta$ that we study there exist at least two groups 
and therefore 
two types of reaction. One of the groups consists of only a single molecule with 
state 0. This particle changes its state from 0 to 1 
spontaneously.  We call this transition an \textit{excitation}. An interpretation of 
this transition can be 
conceived as actions of an external medium, for example, \textit{pumping}. 
Beside, there are other reactions which involve groups of particles. In these reactions  
some amount of excited particles transit simultaneously to 0. Physically, we assume
this transition to be accompanied by a emission of light (radiation). Hence the
interpretation of our model is describing luminescence.%as to is

Our goal is to find statistical properties of the luminescence in
some scaling limits. The first limit is the thermodynamic limit $N\to\infty$ with a 
suitable scaling of the variables. This limit brings us to the area of the large 
deviation theory. We consider the large deviations on the path level providing a 
corresponding \textit{rate function}.  
Having the rate function we can study rare events. Our interest is an event that the cumulative 
number of the transitions from 1 to 0 during a finite time-interval is at least $B$. Then the 
second limit is when $B\to\infty$. We study the density of excited particles conditioned 
that $B$ is very large. In the interpretation as emission of light the large deviation 
describe a bright flash on the time interval. The result we have obtained is rather 
unexpected:
%\vspace{.3cm}
\begin{enumerate}
\item[]
\textit{
The average density of the excited molecules conditioned on very large $B$ does not 
depend on the rates of the transitions of the Markov process.
}
%\vspace{.3cm}
\end{enumerate}

This statement means for example that no matter how small is the rate of the pumping 
(the transitions $0\to 1$), for any values of the  rates of the reactions the density 
of the excited molecules is equal to the same value. For example, for the particular 
case of process $\theta$ published in  \cite{PPSVY} this value is equal $\frac12$.

%Of course, there are cases when this result is wrong. However the co-dimension of 
%these cases is small. \red{\tt Why ``of course''?
%What is co-dimension? Do we really need this remark? I feel we should delete it.
%I find it confusing and it raises the
%question under which conditions the results are wrong.}

The stochastic processes we consider in this article are extremely simple. Our efforts 
are concentrated on the large deviation theory for these processes. More complete 
presentations of the stochastic chemical kinetics theory see in 
\cite{Ar,WH,HJ,H1,H2,BP,F,MPR}. We also mention
an intriguing duality of the pumping to one-dimensional diffusion-limited 
pair-annihilation \cite{Schu95}.

The next section contains all definitions. 

%\newpage

\section{Definitions}

%\subsection{Two Markov Processes}
%\subsubsection*{The Markov process 1}
Here we formally describe the set ${\Theta}$ of the Markov processes. We define $\mathscr N=\{1,...,N\}$. 
The elements of $\mathscr N$ are called particles (or the molecules as in the Introduction). We also define the function
\[
c:\:\mathscr N\to \{0,1\},
\]
where the value $c(i)$ is called a \textit{state} of the particle $i\in\mathscr N$. 
All processes $\theta\in\Theta$ have the same set $\C$ of states which is the set of all configurations  
$\mathbb{C} = \{0,1\}^N$.

The configuration set can be presented as the union 
\bel{3.1a} 
\C=\bigcup_{m=0}^N\C(m),
\ee
where  
 \[
\C(m)=\left\{ c\in\mathbb C:\   \sum_{i=1}^N c(i) =m \right\}.
\]
is a set of all configurations having the same number $m$ of the particles in state 1. 
%The set $\C(m)$ is called $m$th \textit{layer} of the configuration space $\C$.
Every configuration $c$ splits the set $\Ns=\Ns_0(c)\cup\Ns_1(c)$, where 
\bel{3.3a}
\Ns_0(c)=\{i:\:c(i)=0\},\ \Ns_1(c)=\Ns\setminus\Ns_0(c). 
\ee
In other words, for a given configuration $c$, the set $\Ns_\alpha(c)$ is the set of atoms in state $\alpha\in\{0,1\}$.
It is clear that $m=|\Ns_1(c)|$.  

In order to define the transitions of the processes $\theta\in\Theta$
we introduce the triplet
\[
\tau_0 := (1,0,-1)
\]
and the set of integer triplets
\[
\mathcal{K} :=\big\{(k, r, s): \ k, r, s\in\mathbb{N} \mbox{ such that } 
N\ge k\ge r \ge s \big\} \cup \tau_0. 
\]
Observe that in any given triplet the third parameter is either a positive 
integer or $-1$. In the latter case only the triplet $\tau_0$ belongs to $\mathcal{K}$. 
The transitions of a given process $\theta$ are then defined by $d+1$ 
intensities $\mu_0$, $\tilde{\mu}_i$,
$i=1,\dots,d$ and a
finite subset 
\bel{3.3b}
\mathcal {K}_{\theta}=\{\tau_i=(k_i,r_i,s_i) \in \mathcal{K}, i=1,\dots,d:\ 
\mbox{such that }\tau_i\ne\tau_j,\mbox{ if }i\ne j\} \cup \{\tau_0\} \subset \mathcal{K}
\ee
of $d+1$ distinct triples from $\mathcal{K}$ where each triplet encodes an allowed
transition as follows.

\begin{itemize}
\item Pumping: Consider a configuration $c\in\C(m)$ with $0\le m<N$, which
means that $N-m$ particles are in the ground state. Then the transition corresponding to 
the triplet $\tau_0=(1,0,-1)$ means that one of the particles in the ground state, which is
chosen uniformly, becomes 
excited. In our notation the first entry 1 in the triplet means that one particle takes part 
in the transition, the second 0 means that this particle is in the ground state before the 
transition takes place, and -1 denotes that this particle is being excited. Physically, this transition is due to external pumping and occurs with intensity $\mu_0$.

\item Radiation: The triplets $\tau_i$ with $i\in\{1,\dots,d\}$
denote transitions of a configuration $c\in\C(m)$ 
with $0<m\leq N$ where the number of the excited particles decreases as follows. 
The first parameter $k_i$ is the {\it cluster size}, i.e., the number of the particles 
involved in the transition $\tau_i$. A cluster is chosen uniformly from all possible
clusters of $k_i \leq N$ particles. If $r_i$ of the particles in the selected cluster
are excited then $s_i$ uniformly chosen particles from this set relax to their ground state. 
Otherwise nothing happens. We call this event radiation and denote its intensity
by $\tilde{\mu}_i$. Physically, this process is a spontaneous emission of $s_i$ photons
that occurs as a result of a collective interaction of $k$ atoms.
\end{itemize}

Notice that these processes induce transitions between the sets $\C(m)$.
Pumping increases the number of excited by one.
Since there are $N-m$ ground state particles that can turn excited 
the intensity of the transition from $c\in\C(m)$ to the set $\C(m+1)$ is equal to 
$\mu_0(N-m)$. Radiation defined by some triplet $\tau_i$ moves the configuration 
$c$ from the set $\C(m)$ to the set $\C(m-s_i)$. Since
for a fixed $m$ the number clusters with $k_i$ particles of which $r_i$ 
are excited is given by
\[
q_{i}(N,m) = \binom{m}{r_i}\binom{N-m}{k_i-r_i}
\]
the intensity of the transition $\C(m)$ to $\C(m-s_i)$ is $\tilde{\mu}_i q_{i}(N,m)$.
Below, the rates $\tilde{\mu}_i$ will be chosen to depend on $N$, 
$\tilde{\mu}_i\equiv \tilde{\mu}_i(N)$. 
We define this functional dependence when we determine the Markov process of 
functionals of the process $\theta$.
  
Let us introduce some auxiliary objects corresponding to these $d$ transitions. 
Let $S,R,K$ be the following three subsets associated to the configuration $c$:  
\begin{enumerate}
\item[]
$S,R\subseteq \Ns_1(c), K\subseteq \Ns$ such that $\varnothing \ne S\subseteq R\subseteq K$
\end{enumerate}
We call these sets \textit{three-set-$c$} and denote as $\langle S,R,K\rangle^c$ or by the shorter notation $\pmb t^c=\langle S,R,K\rangle^c$. We use also the equivalent presentation $\pmb t^c=\langle\mathbf s(\pmb t^c),\mathbf r(\pmb t^c),\mathbf k(\pmb t^c)\rangle^c$ with the obvious notations.

We say that the three-set-$c$ $\pmb t^c$ represents the triplet $\tau_i=(k_i,r_i,s_i)\in\mathcal K_\theta$ and write $\pmb t^c\to\tau_i$ if $|S|=s_i,|R|=r_i,|K|=k_i$. The index $i$ in this representation we denote as $i(\pmb t^c)$. Let $\mathbf s(\pmb t^c)=S$.
Every three-set-$c$ $\pmb t^c=\langle S,R,K\rangle^c$ representing the reaction $\tau_i=(k_i,r_i,s_i)$ realizes the transition $c\to c'=c_{S_i}^-$, where 
$$
c_{S_i}^-(j)=\begin{cases}0,&\mbox{ if }j\in S_i\\c(j),&\mbox{ otherwise}\end{cases}.
$$
On the other hand, the process represented by the triplet $\tau_0$ realizes the transition
$c\to c'=c_{j_0}^+(j)$, where 
$$
c_{j_0}^+(j)=\begin{cases}1,&\mbox{ if }j=j_0\\c(j),&\mbox{ otherwise}\end{cases}.
$$
Thus, the infinitesimal generator of the process $\theta$ described informally above
is
\[
\mathbf{L}_\theta f(c)=\mu_0\sum_{j_0\in\Ns_0(c)}[f(c_{j_0}^+)-f(c)]+
\sum_{S\subseteq\Ns_1(c)}\sum_{\tau\in\mathcal K_\theta}\;\sum_{\substack{\pmb t^c:\:\pmb t^c\to\tau,\\
\mathbf s(\pmb t^c)=S}}\mu_{i(\pmb t^c)}[f(c_{S}^-)-f(c)] .
\]

\textbf{Functionals}.
Introduce some functionals of the processes $\theta$ we deal with in this work. The first functional is simply defined as the number of $1$'s in a current configuration
$$
m(t) = \bigl|\,  \Ns_1 \bigl( \theta(t) \bigr)\,  \bigr|.
$$

We consider the processes on time interval $[0,T]$.
Let $0\equiv t_0 < t_1 < ... < t_n < T$ be the sequence of the moments on the time interval $[0,T]$ where the process $\theta$ change its configuration. That is the sequence $\theta(t_0)=c_0,  \theta(t_1)=c_1,...,\theta(t_n)=c_n,$ is the path of the $\theta$ on $[0,T]$: assuming $\theta(t) = c_i, t_i\le t < t_{i+1}$, where $t_{n+1} \equiv T$, i.e. the process $\theta$ is continuous from the right, therefore $\theta(t_k)$ is the configuration $c_k$ which the process takes after the random event occurred at $t_k$. 

Any jump-moment $t_k, k=1,\dots, n$ corresponds to some of the triplets $\tau\in\mathcal K_\theta$. For any $i=1,\dots, d$ introduce $I^i(t)$ the number of jumps of the process  $\theta$ occurred thanks the triplet $\tau_i$. Define the following $d$ functionals:
\[
y^i(t)=s_i I^i(t),\ i=1,\dots, d.
\]
All $d+1$ functionals we deal with can be considered as a vector path $(m(t),y^1(t),...,y^d(t))$.

Note that the vector of functional $(m(t),y^1(t),...,y^d(t))$ is indeed Markov process on the state space $\{0, 1, \dots, N\} \times \mathbb{Z}_+^d$. We denote it as $\zeta=(\xi,\eta_1,...,\eta_d)$.
 The continuous time Markov process $\zeta$ is defined by a given set of triples $\mathcal{K}_\theta$ and its dynamics  governed by the generator 
\begin{equation}\label{gen}
\begin{aligned}
\mathbf{L}g(m,y_1,...,y_d)&=\mu_0(N-m)[g(m+1,y_1,...,y_d)-g(m,y_1,...,y_d)]\\
&+\sum_{i=1}^d \tilde\mu_i \binom{m}{r_i}\binom{N-m}{k_i-r_i}  \bigl[ g(m-s_i,y_1,...,y_i+s_i,...,y_d)-g(m,y_1,...,y_d) \bigr],
\end{aligned}
\end{equation}
acting on the set of functions $g:\;\N\times \Z_+^d \to\R$, where the parameter $\mu_0$ is the rates of the jumps $+1$ of the process component $\xi$, and $\tilde\mu_i,i=1,...,d$ are rates of the jumps $-s_i$ of the process $\xi$. When $\xi$ changes its value by $-s_i$, simultaneously the component $\eta_i$ increases its value by $s_i$: every component $\eta_i$ jumps by the non-random number $s_i$, defined by the triplet $(k_i,r_i,s_i)$. 

Except for the case $\mathcal K_\theta=\{\tau_1=(1,1,1), \tau_0=(1,0,-1) \}$ (the linear case) all processes describe the mean field behavior. Thus the intensity of  the emissions $\tilde\mu_i \equiv \tilde\mu_i(N)$ are scaled in the limit $N\to\infty$
\begin{equation}\label{3.1b}
\mu_i=\lim_{N\to\infty} \tilde\mu_i(N) N^{k_i-1},
\end{equation}
where $\mu_i$ is some positive constant, or we can simply define
\bel{3.1.1}
\tilde\mu_i \equiv \tilde\mu_i (N) = \frac{\mu_i}{N^{k_i-1}}.
\ee
This type of mean field scaling sometimes is called \textit{ canonical } scaling, see for example \cite{MPR}, or it refers to the stochastic analog of the \textit{ law of mass actions}, see for example \cite{EK}, Chapter 10, pp. 454.  Moreover, such rate scaling brings our process
defined by the generator (\ref{gen}) to the class of the processes known as \textit{density dependent processes}, see \cite{EK}, Chapter 10. Indeed, with \reff{3.1.1} the generator (\ref{gen}) can be represented in the following way
\begin{equation}\label{gen1}
\begin{aligned}
&\mathbf{L}g(m,y_1,...,y_d)=N\mu_0 \Bigl(1-\frac{m}{N} \Bigr)[g(m+1,y_1,...,y_d)-g(m,y_1,...,y_d)]\\
&+\sum_{i=1}^d N \left( \mu_i \frac{(m/N)^{r_i}}{r_i!} \frac{(1-m/N)^{k_i-r_i}}{(k_i-r_i)!}   + O\Bigl(\frac1N\Bigl) \right) \bigl[ g(m-s_i,y_1,...,y_i+s_i,...,y_d)-g(m,y_1,...,y_d) \bigr].
\end{aligned}
\end{equation}
Thus, our process will belong to the \textit{density dependent family} of processes introduced \cite{EK}, if for any $i=1, \dots, d$ we define $\tilde\mu_i(N)$ such that
\bel{rate.scaling}
\tilde\mu_i(N) \binom{m}{r_i}\binom{N-m}{k_i-r_i} =
N \mu_i \frac{(m/N)^{r_i}}{r_i!} \frac{(1-m/N)^{k_i-r_i}}{(k_i-r_i)!}
\ee
for any $N$. This scaling avoids $O(1/N)$ in (\ref{gen1}). It provides finally the generator for the process $\zeta$ we will deal with in this paper:
\begin{equation}\label{gen2}
\begin{aligned}
&\mathbf{L}g(m,y_1,...,y_d)=N\mu_0 \Bigl(1-\frac{m}{N} \Bigr)[g(m+1,y_1,...,y_d)-g(m,y_1,...,y_d)]\\
&+ N\sum_{i=1}^d \mu_i \frac{(m/N)^{r_i}}{r_i!} \frac{(1-m/N)^{k_i-r_i}}{(k_i-r_i)!} \bigl[ g(m-s_i,y_1,...,y_i+s_i,...,y_d)-g(m,y_1,...,y_d) \bigr].
\end{aligned}
\end{equation}
Further, the $\mathbf{L}$ stands for the generator defined by (\ref{gen2}).

\section{Problems and the approach}

We  study the large deviations  as the following. Namely, we want to find the asymptotics of
\[
\ln\Pr\left(\sum_{i=1}^d\eta_i(T)>BN\right)
\]
as $N\to\infty$, where $B>0$. 
Or, equivalently, 
\begin{equation}\label{2.1}
\ln\Pr\left(\frac1N\sum_{i=1}^d\eta_i(T)>B\right).
\end{equation} 

Since the problem belongs to the branch of the probability theory called the large deviations we have to find the rate function corresponding to the large deviation principle related to the problem (\cite{DZ}, \cite{FeK}, \cite{Reza}, see also \cite{K}, where the large deviation principle is studied for a model that is close to our model in some respect).

The rate function $I(\cdot)$ is defined on the space of paths $\gamma(\cdot) = (x_0(\cdot),x_1(\cdot),...,x_d(\cdot))$ mapping $[0,T]\to[0,1]\times \R_+^d$, where the paths $x_1(\cdot),...,x_d(\cdot)$ start from 0 at $t=0$ and are non-decreasing such that:
 $$
\Pr \left( \frac{ \zeta (\cdot)}{N} \approx \gamma(\cdot) \right) \approx \exp \left( - N I\bigl( \gamma(\cdot)\bigr) \right). % (\int_0^T L(\gamma(t), \dot{\gamma}(t) ) dt \right)
$$
We define the rate function via non-linear Hamiltonian
\begin{equation}\label{3.3c}
%\begin{array}{lll}
\begin{aligned}
H_{\mathcal K_\theta} & \bigl(x_0(\cdot),x_1(\cdot),...,x_d(\cdot),\sigma(\cdot),\k_1(\cdot),...,\k_d(\cdot) \bigr) = \\
&{} = \mu_0(1-x_0(\cdot))[e^{\sigma(\cdot)}-1]+\sum_{i=1}^d \mu_i\frac{x_0^{r_i}}{r_i!}\frac{(1-x_0)^{k_i-r_i}}{(k_i-r_i)!}[e^{-s_i\sigma(\cdot)+s_i\k_i(\cdot)}-1].
\end{aligned}
%\end{array}
\end{equation}

The Hamiltonian $H_{\mathcal K_\theta}$ is defined by the generator $\mathcal{H}$ of a non-linear semigroup acting on the functions $f$ on $[0,1]\times \R_+^d$. Let $\underline{x} \in [0,1]\times \R_+^d$ then
$$
H_{\mathcal K_\theta} (\underline{x}, \triangledown\! f) := (\mathcal{H} f)(\underline{x}) :=\lim_{N\to\infty} \frac 1 N e^{-N f(\underline{x})} \left(  \mathbf{L} e^{Nf} \right) (\underline{x})
$$
The rate function is defined as the integral of the Legendre transform of the Hamiltonian.
\begin{eqnarray}\label{2.2}
I\left(x_0(\cdot), x_1(\cdot),...,x_d(\cdot)\right) 
&=& \int_0^T\sup_{\sigma(t),\k_1(t),...,\k_d(t)}\Bigl[ \sigma(t)\dot x_0(t)+\sum_{i=1}^d\k_i(t)\dot x_i(t)  \nonumber \\
&& {} -  H_{\mathcal K_\theta} \bigl( x_0(t),x_1(t),...,y_d(t),\sigma(t),\k_1(t),...,\k_d(t) \bigr) \Bigr] \ed t
\end{eqnarray}
The asymptotics of the probability in \reff{2.1} is
\[
\lim_{N\to\infty} \frac 1 N \ln\Pr\left(\frac1N\sum_{i=1}^d\eta_i(T)\geq B\right)=-\inf \left\{I\left(x_0,x_1,...,x_d\right):\:\sum_{i=1}^d x_i(T)\geq B\right\}.
\]
Finding this asymptotics we have to solve the following system
\begin{equation}\label{hameq}
\left\{ 
\begin{aligned}
\dot x_0=&\mu_0 (1-x_0)\exp\{\sigma\} - \sum_{i=1}^d \mu_i s_iQ(x_0,r_i)Q(1-x_0,k_i-r_i) \exp\{-(s_i\sigma-s_i\varkappa_i)\} , \\
\dot x_1=&\mu_1 s_1\exp\{-(s_1\sigma-s_1\varkappa_1)\}, \\
&...\\
\dot x_d=&\mu_d s_d\exp\{-(s_d\sigma-s_d\varkappa_d)\}, \\
\dot \sigma=&\mu_0 [\exp\{\sigma\}-1]-  \sum_{i=1}^d \mu_i [\exp\{-(s_i\sigma-s_i\varkappa_i)\} -1] \times  \\
& {} \times [Q(x_0,r_i-1)Q(1-x_0,k_i-r_i)-Q(x_0,r_i)Q(1-x_0,k_i-r_i-1)] \\
\dot\varkappa_1=&0\\
&...\\
\dot\varkappa_d=&0,
\end{aligned} \right.
%\end{array} \right.
\end{equation}
where 
$Q(z,q)=\frac{z^{q}}{q!}$, if  $z\in[0,1],q\in\Z_+$ and $Q(z,q)=0$ elsewhere.  The first $d+1$ equations of (\ref{hameq}) are to find the supremum in the integrand (\ref{2.2}). The remaining $d+1$ equations are the system of Euler-Lagrange equations. 

The explicit solution of the Hamiltonian system \reff{hameq} is either very difficult or impossible for the case general initial conditions $x_0(0)$. The problem appears from the expression $Q(x_0,r_i-1)Q(1-x_0,k_i-r_i)$ and $Q(x_0,r_i)Q(1-x_0,k_i-r_i-1)$ in the equation for $\sigma$. There exists perhaps only exclusion for the case studied in \cite{PPSVY}, where the exact solution is found (see also the examples below).

However there is an initial conditions $x_0(0)=\wh{x}_0$ such that the solution can be found. This solution is such that $x_0(t)\equiv \wh x_0,\ t\in[0,T]$. 

The main result in this note is about the existence of such solution.
The main goal of our studies is to understand what on level $B\to\infty$ happens. A general conjecture is
\begin{theorem}
There exists only one $i_0\in\{1,...,d\}$ such that 
\[
x_0(t)\equiv \wh{x}_0=\frac{r_{i_0}}{k_{i_0}+s_{i_0}}
\]
for $0\leq t\leq T$.
\end{theorem}

Of course, some assumption are required.  Later only particular cases are considered.

%\newpage

\section{Proof} We split the proof in several steps corresponding different values of $d$.
\subsection{$d=1$} We prove the theorem for the processes with one reaction.
Let $\mathcal{K}_\theta=\{(k,r,s), (1,0,-1)\}$. The Hamiltonian is
\bel{hamd1}
H(x_0,x_1,\sigma,\k)=\mu_0(1-x_0)[e^\sigma-1]+Q(x_0,r)Q(1-x_0,k-r)[e^{-s\sigma+s\k}-1],
\ee
%Further we consider two cases: 1) $k-r\geq 1$, 2)$k-r=0$. 
The Hamiltonian system (\ref{hameq}) for this case is
\begin{equation}\label{8.2}
\left\{
\begin{array}{lll}
\dot x_0&=&\mu_0 (1-x_0)e^\sigma-sQ(x_0,r)Q(1-x_0,k-r)e^{-s\sigma+s\k}\\
\dot x_1&=&sQ(x_0,r)Q(1-x_0,k-r)e^{-s\sigma+s\k}\\
\dot\sigma&=&\mu_0(e^\sigma-1)-\bigl(Q(x_0,r-1)Q(1-x_0,k-r)-Q(x_0,r)Q(1-x_0,k-r-1)\bigr)[e^{-s\sigma+s\k}-1]\\
\dot\k&=&0
\end{array} \right.
\end{equation}
Recall that $Q(1-x_0,k-r-1)=0$ in the case $k=r$.

Assume that there exists $\wh{x}_0$ such that $x_0(t)\equiv \wh{x}_0$,  is a constant. Assume that $\dot x_1$ is also a constant. Then $\dot x_1=B$.
Then, assuming  $T=1$, we obtain %\red{(why $\dot{y} = B$?)}
\[
e^{-s\sigma+s\k}=\frac B{sQ(\wh{x}_0,r)Q(1-\wh{x}_0,k-r)}.
\]
Since $\dot x_0=0$ we obtain
\[
\mu_0 e^\sigma=\frac{sQ(\wh{x}_0,r)Q(1-\wh{x}_0,k-r)e^{-s\sigma+s\k}}{1-\wh{x}_0}=\frac B{1-\wh{x}_0}
\]
Because $\dot\sigma=0$ we have
\[\mu_0 e^\sigma=
\mu_0+\left(Q(\wh{x}_0,r-1)Q(1-\wh{x}_0,k-r)-Q(\wh{x}_0,r)Q(1-\wh{x}_0,k-r-1)\right)[e^{-s\sigma+s\k}-1]
\]
From three above relations we obtain
\[
\frac1{1-\wh{x}_0}=\frac{Q(\wh{x}_0,r-1)Q(1-\wh{x}_0,k-r)}{sQ(\wh{x}_0,r)Q(1-\wh{x}_0,k-r)}-\frac{Q(\wh{x}_0,r)Q(1-\wh{x}_0,k-r-1)}{sQ(\wh{x}_0,r)Q(1-\wh{x}_0,k-r)}+o(1)
\]
for $B\to\infty$. This is equivalent to
\[
1=\frac rs\frac{1-\wh{x}_0}{\wh{x}_0}-\frac{k-r}s
\]
which finally gives
\bel{4.1}
x_0(t) \equiv \wh{x}_0 =\frac r{k+s}.
\ee

%\subsection{Examples}
We consider several examples for the demonstration of the result.
%In further we drop the tilda over $\mu$ (see \reff{3.1})

\subsubsection{Linear case}
It is the case when $\mathcal K_\theta=\{ (1,1,1), (1,0,-1)\}$. The Hamiltonian is
\[
%\begin{array}{lll}
H_{\mathcal K_\theta}(x_0(\cdot),x_1(\cdot),\sigma(\cdot),\k(\cdot))=\mu_0 (1-x_0(\cdot))[e^{\sigma(\cdot)}-1]+\mu_1x_0(\cdot)[e^{-\sigma(\cdot)+\k(\cdot)}-1].
%\end{array}
\]
The corresponding Hamiltonian system (\ref{8.2}) is
\begin{equation}\label{3.1c}
\left\{ \begin{array}{rcl}
\dot x_0&=&\mu_0 (1-x_0)\exp\{\sigma\}+ \mu_1 x_0\exp\{ -\sigma+\k\}, \\
\dot x_1&=&\mu_1\exp\{-(\sigma-\varkappa)\}, \\
\dot \sigma&=&\mu_0 [\exp\{\sigma\}-1]- \mu_1 [\exp\{-(\sigma-\varkappa)\} -1], \\
\dot\varkappa&=&0\\
\end{array} \right.
\end{equation}
This system can be solved analytically for the following initial and boundary conditions
\bel{inibc}
\begin{array}{lcl}
x_0(0)&=& \wh{x}_0 \in[0,1]\\
x_1(0)&=&0\\
x_1(T)&=&B\\
\sigma(T)&=&0.
\end{array}
\ee
There is no the mean field actions in this case therefore the parameter $\tilde\mu_1=\mu_i$ is not scaled. According \reff{4.1}
\[
x_0(t)=\frac12.
\]
This case has been studied in \cite{PPSVY}.

%\subsection{Non-linear cases} 
\subsubsection{Quadratic potentials}\label{quadrat} In this section we consider three examples: $\tau_1 = (2,2,2), \tau_1=(221)$ and $\tau_1 = (2,1,1)$. 
%The absorption in all three examples is described in the same way as it  was in the general cases. The difference of the processes are in the actions of the emission.

The Hamiltonian equations \reff{hameq} of these examples cannot be solved explicitly for any initial condition $x_0(0)$ (as we assume). However there exists a special $x_0(0)=\wh x_0$ such that the solution $x_0(t)$ is a constant, $x_0(t)\equiv\wh x_0$ if $x_0(0)=\wh x_0$. We  obtain the value of $\wh x_0$   in the course of the solution.

The following Hamiltonians correspond to the processes considered here.
\begin{enumerate}
\item 
$\mathcal K_\theta=\{(222), (1,0,-1)\}$. The corresponding Hamiltonian \reff{hamd1} in this case is

$$
H_{\mathcal K_\theta}(x_0,x_1,\sigma,\k)=\mu_0(1-x_0)[e^{\sigma}-1]+\mu_1\frac{x_0^2}{2}[e^{-2\sigma+2\k}-1].
$$

The interpretation of this process emission is the following. The emission occurs when two excited atoms collide each other. Each collided atom emits one photon. Therefore the second term in the Hamiltonian has multiplier $x_0^2/2$. By the same reason, the exponent is $e^{-2\sigma+2\k}$ meaning that the number of the excited atoms decreased by two and two photons were emitted. The corresponding Hamiltonian system \reff{8.2} is
$$
\left\{
\begin{array}{lll}
\dot x_0&=&\mu_0(1-x_0)e^{\sigma}-\mu_1 \frac{x_0^2}{2}e^{-2\sigma+2\k}\\
\dot x_1&=&\mu_1 x_0^2e^{-2\sigma+2\k}\\
\dot\sigma&=&\mu_0[e^{\sigma}-1]-\mu_1 x_0[e^{-2\sigma+2\k}-1]\\
\dot\k&=&0
\end{array}
\right.
$$
with the same boundary conditions \reff{inibc}, and according  \reff{4.1}
%\[
%\begin{array}{lll}
%x_0(0)=\wh x_0,\\
%x_1(0)=0,\\
%x_1(T)=B,\\
%\sigma(T)=0.
%\end{array}
%\]
%We assume that $x^0$ is such that $x_1(t)\equiv x^0$, i.e is the constant on $[0,T]$.
in this case 
\[
x_0(t)\equiv \wh x_0 = \frac12.
\]

\item $\mathcal K_\theta=\{(2,2,1), (1,0,-1)\}$. The corresponding Hamiltonian \reff{hamd1} in this case is
$$
H_{\mathcal K_\theta}(x_0,x_1,\sigma,\k)=\mu_0(1-x_0)[e^{\sigma}-1]+\mu_1\frac{x_0^2}{2}[e^{-\sigma+\k}-1].
$$
As in the previous process, the emission occurs when two excited atoms collide each other. However only one of the collided atoms emits photon. Therefore there is  the multiplier $x_0^2/2$ as in the previous case. In this case, the exponent is $e^{-\sigma+\k}$ meaning that the number of the excited atoms is decreasing by one and one photon is emitted.
The Hamiltonian system \reff{8.2} is
$$%\begin{equation}\label{5.1}
\left\{
\begin{array}{lcl}
\dot x_0&=&\mu_0 (1-x_0)e^{\sigma}-\mu_1 \frac{x_0^2}{2}e^{-\sigma+\k},\\
\dot x_1&=&\mu_1 \frac{x_0^2}{2}e^{-\sigma+\k},\\
\dot\sigma&=&\mu_0 [e^{\sigma}-1]-\mu_1 x_0[e^{-\sigma+\k}-1],\\
\dot\k&=&0
\end{array}
\right.
$$%\end{equation}
with the same boundary conditions \reff{inibc}, but
%\[
%\begin{array}{lll}
%x_0(0)=\wh x_0,\\
%x_1(0)=0,\\
%x_1(T)=B,\\
%\sigma(T)=0.
%\end{array}
%\]
in this case 
\[
x_0(t)\equiv \wh x_0 =\frac23.
\]

\item $\mathcal K_\theta=\{(2,1,1), (1,0,-1)\}$. The corresponding Hamiltonian \reff{hamd1} in this case is
$$
H_{\mathcal K_\theta}(x_0,x_1,\sigma,\k)=\mu_0(1-x_0)[e^{\sigma}-1]+\mu_1 x_0(1-x_0) [e^{-\sigma+\k}-1].
$$

In this case the reaction due to two atoms, one of them is excited and another one is in the ground state. Therefore the multiplier is $x_0(1-x_0)$. The exponent is the same as in the previous case: $e^{-\sigma+\k}$. The Hamiltonian system is
\[
\left\{
\begin{array}{lcl}
\dot x_0&=&\mu_0(1-x_0)e^{\sigma}-\mu_1 x_0(1-x_0)e^{-\sigma+\k},\\
\dot x_1&=&\mu_1 x_0(1-x_0)e^{-\sigma+\k},\\
\dot\sigma&=&\mu_0 [e^{\sigma}-1]-\mu_1 (2x_0-1)[e^{-\sigma+\k}-1],\\
\dot\k&=&0
\end{array}
\right.
\]
with the same boundary conditions \reff{inibc} and 
%\[
%\begin{array}{lll}
%x_0(0)=\wh x_0,\\
%x_1(0)=0,\\
%x_1(T)=B,\\
%\sigma(T)=0.
%\end{array}
%\]
in this case 
\[
x_0(t)=\wh x_0 = \frac13.
\]

\end{enumerate}

\subsubsection{Cubic potentials} In this section we consider a model having cubic potential which means that the emission occurs when three atoms collide.
\begin{enumerate} 
\item[] If $\mathcal K_\theta=\{(333),(1,0,-1)\}$, then
$$
H_{\mathcal K_\theta } (x_0,x_1,\sigma,\k)=\mu_0 (1-x_0)[e^{\sigma}-1]+\mu_1 \frac{x_0^3}{6}[e^{-3\sigma+3\k}-1].
$$
Three photons are emitted if three excited atoms collide. The number of possibilities to choose three excited atoms gives the multiplier $x_0^3/6$ and the exponent is $e^{-3\sigma+3\k}$ meaning that three excited atoms  disappear  and three photons were emitted. The Hamiltonian system is
\[
\left\{
\begin{array}{lcl}
\dot x_0&=&\mu_0(1-x_0)e^{\sigma}-\mu_1 \frac{x_0^3}{2}e^{-3\sigma+3\k},\\
\dot x_1&=&\mu_1 \frac{x_0^3}{2}e^{-3\sigma+3\k},\\
\dot\sigma&=&\mu_0 [e^{\sigma}-1]-\mu_1 \frac{x_0^2}{2}[e^{-3\sigma+3\k}-1],\\
\dot\k_2&=&0
\end{array}
\right.
\]
with the boundary conditions \reff{inibc} and  
%\[
%\begin{array}{lll}
%x_0(0)=\wh x_0,\\
%x_1(0)=0,\\
%x_1(T)=B,\\
%\sigma(T)=0.
%\end{array}
%\]
%In this case 
\[
x_0(t)\equiv \wh x_0 =\frac12.
\]

\end{enumerate}

\subsection{$d=2$} 

In this section we consider two reactions. 
Let
\[
\mathcal K_\theta=\{(k_1,r_1,s_1), (k_2,r_2,s_2), (1,0,-1)\}
\]
and $s_1\ne s_2$. The Hamiltonian is
 
 \begin{equation}\label{7.0}
 \begin{array}{lcl}
 H_{\mathcal K_\theta}(x_0,x_1,x_2,\sigma,\k_1,\k_2)&=&
 \mu_0(1-x_0)[e^{\sigma}-1]+\mu_1Q(x_0,r_1)Q(1-x_0,k_1-r_1)[e^{-s_1\sigma+s_1\k_1}-1]
 \\
&+& \mu_2Q(x_0,r_2)Q(1-x_0,k_2-r_2)[e^{-s_2\sigma+s_2\k_2}-1]
\end{array}
\end{equation}
The Hamiltonian system is
\begin{equation}\label{7.1}
\left\{
\begin{array}{lll}
\dot x_0&=&\mu_0(1-x_0)e^{\sigma} \\ 
&& {} -\mu_1s_1 Q(x_0,r_1)Q(1-x_0,k_1-r_1)e^{-s_1\sigma+s_1\k_1} \\
&& {} -\mu_2s_2Q(x_0,r_2)Q(1-x_0,k_2-r_2)e^{-s_2\sigma+s_2\k_2},\\
\dot x_1&=&\mu_1s_1Q(x_0,r_1)Q(1-x_0,k_1-r_1)e^{-s_1\sigma+s_1\k_1},\\
\dot x_2&=&\mu_2s_2Q(x_0,r_2)Q(1-x_0,k_2-r_2)e^{-s_2\sigma+s_2\k_2},\\
\dot\sigma&=&\mu_0[e^{\sigma}-1] \\
&& {} -\mu_1\left[Q(x_0,r_1-1)Q(1-x_0,k_1-r_1)-Q(x_0,r_1)Q(1-x_0,k_1-r_1-1)\right][e^{-s_1\sigma+s_1\k_1}-1]\\
&&{} -\mu_2\left[Q(x_0,r_2-1)Q(1-x_0,k_2-r_2)-Q(x_0,r_2)Q(1-x_0,k_2-r_2-1)\right][e^{-s_2\sigma+s_2\k_2}-1],\\
\dot\k_1&=&0,\\
\dot\k_2&=&0.
\end{array}
\right.
\end{equation}
The boundary conditions are
\[
\begin{array}{lcl}
x_0(0)&=&\wh x_0,\\
x_1(0)&=&0,\\
x_1(T)&=&B_1,\\
x_2(0)&=&0,\\
x_2(T)&=&B_2,\\
\sigma(T)&=&0,
\end{array}
\]
where $B=B_1+B_2$.

Let $i_0=\arg(\max\{s_1,s_2\})$. We prove that
\begin{equation}\label{8.0}
x_0(t)=\frac{r_{i_0}}{k_{i_0}+s_{i_0}}.
\end{equation}

As in previous case we assume that $x_0(t)\equiv \wh x_0$. Then $\sigma\equiv constant$, and $\k_1\equiv constant$ and $\k_2\equiv constant$. And
\begin{equation}\label{8.1}
\begin{aligned}
\dot x_1&=B_1=\mu_1s_1Q(x_0,r_1)Q(1-x_0,k_1-r_1)e^{-s_1\sigma+s_1\k_1}\\
\dot x_2&=B_2=\mu_2s_2Q(x_0,r_2)Q(1-x_0,k_2-r_2)e^{-s_2\sigma+s_2\k_2}
\end{aligned}
\end{equation}

Recall that (see the first equation in \reff{7.1})
\begin{equation}\label{8.2.1}
\begin{aligned}
0=\dot x_0=\mu_0(1-x_0)e^{\sigma}
&-\mu_1s_1Q(x_0,r_1)Q(1-x_0,k_1-r_1)e^{-s_1\sigma+s_1\k_1} \\
&-\mu_2s_2Q(x_0,r_2)Q(1-x_0,k_2-r_2)e^{-s_2\sigma+s_2\k_2}
\end{aligned}
\end{equation}
and (see the fourth equation in \reff{7.1})
\begin{equation}\label{8.3}
\begin{aligned}
0=\dot\sigma=&\mu_0[e^\sigma-1] \\
&-\mu_1\left[Q(x_0,r_1-1)Q(1-x_0,k_1-r_1)-Q(x_0,r_1)Q(1-x_0,k_1-r_1-1)\right][e^{-s_1\sigma+s_1\k_1}-1]\\
&-\mu_2\left[Q(x_0,r_2-1)Q(1-x_0,k_2-r_2)-Q(x_0,r_2)Q(1-x_0,k_2-r_2-1)\right][e^{-s_2\sigma+s_2\k_2}-1]
\end{aligned}
\end{equation}

%$\cfrac{T}{M}$

Finding relation between $B_1$ and $B_2$ we assume that $B_1=\alpha B$ and $B_2=(1-\alpha)B$, where $\alpha\in[0.1]$.

Using \reff{8.1} we obtain  from \reff{8.2.1} the  relation
\[\begin{array}{lll}
\cfrac{\mu_0 e^\sigma}{B}&=&
\cfrac{\alpha}{s_1}\left[\cfrac{Q(x_0,r_1-1)Q(1-x_0,k_1-r_1)}{Q(x_0,r_1)Q(1-x_0,k_1-r_1)}-\cfrac{Q(x_0,r_1)Q(1-x_0,k_1-r_1-1)}{Q(x_0,r_1)Q(1-x_0,k_1-r_1)}\right] \\[.5cm]
&&+ \cfrac{1-\alpha}{s_2}\left[\cfrac{Q(x_0,r_2-1)Q(1-x_0,k_2-r_2)}{Q(x_0,r_2)Q(1-x_0,k_2-r_2)}-\cfrac{Q(x_0,r_2)Q(1-x_0,k_2-r_1-2)}{Q(x_0,r_2)Q(1-x_0,k_2-r_2)}\right]+o(1)
\end{array}
\]
as $B\to\infty$. This is
\[
\frac{\mu_0 e^\sigma}{B}=\frac{\alpha}{s_1}\left[\frac{r_1}{x_0}-\frac{k_1-r_1}{1-x_0}\right]+\frac{1-\alpha}{s_2}\left[\frac{r_2}{x_0}-\frac{k_2-r_2}{1-x_0}\right]+o(1).
\]
From \reff{8.2.1} we obtain
\begin{equation}\label{9.0}
\frac{\mu_0 e^\sigma}{B}=\frac\alpha{1-x_0}+\frac{1-\alpha}{1-x_0}
\end{equation}
Finally
\[
\frac\alpha{s_1}\left[{r_1}\frac{1-x_0}{x_0}-(k_1-r_1)\right]+\frac\alpha{s_1}\left[{r_1}\frac{1-x_0}{x_0}-(k_1-r_1)\right]=1
\]
 Then
\begin{equation}\label{9.1}
\frac{1-x_0}{x_0}=\frac{1+\frac\alpha{s_1}(k_1-r_1)+\frac{(1-\alpha)}{s_2}(k_2-r_2)}{\frac\alpha{s_1} r_1+\frac{(1-\alpha)}{s_2}r_2}.
\end{equation}
Finding $\alpha$ we minimise the rate function $I(x_0,x_1,x_2)$ over $\alpha$ (see \reff{2.2}).  From \reff{9.0} and \reff{8.1} we obtain
\[
\begin{array}{lll}
\sigma&=&\ln B+O(1)\\
\k_1&=&\frac1{s_1}\ln B_1+O(1)=\left(1+\frac1{s_1}\right)\ln B+O(1)\\
\k_2&=&\frac1{s_2}\ln B_2+O(1)=\left(1+\frac1{s_2}\right)\ln B+O(1)
\end{array}
\]
The functions $\dot x_0, \dot x_1, \dot x_2$ and $\sigma,\k_1,\k_2$ are constants. Therefore the rate function can be represented as
\begin{equation}\label{9.3}
I(x_0,x_1,x_2)=\sup_{\sigma,\k_1,\k_2}\{\dot x_1\k_1+ \dot x_2\k_2-H_{\mathcal K_{\theta}}(x_0, x_1, x_2,\sigma,\k_1,\k_2)\}.
\end{equation}
Recall that  $\dot x_0=0$.  Thus
\begin{equation}\label{9.2}
\frac{I(x_0,x_1,x_2)}{B\ln B}=\alpha\left(1+\frac1{s_1}\right)+(1-\alpha)\left(1+\frac1{s_2}\right)+o(1)
\end{equation}
since 
$$
\cfrac{H_{\mathcal K_{\theta}}(x_0, x_1, x_2,\sigma,\k_1,\k_2)}{B}=o(\ln B).
$$

The result \reff{8.2.1} follows from \reff{9.2}.

%\newpage

\section*{Acknowledgement}

The research of E. Pechersky, S. Pirogov and A. Vladimirov was carried out
at the Institute for Information Transmission Problems (IITP), Russian
Academy of Science and funded by the Russian Science Foundation (grant
14-50-00150). 

G.M. Sch\"utz  thanks the IME at the University of S\~ao Paulo for kind hospitality.
His contribution was financed in part by Coordenação de Aperfeiçoamento de 
Pessoal de Nível
Superior -- Brazil (CAPES) -- Finance Code 001, and by the grants 
2017/20696-0, 2017/10555-0 of the S\~ao Paulo Research Foundation (FAPESP).

A. Yambartsev thanks  Conselho Nacional de Desenvolvimento
Cient\'ifico e Tecnol\^ogico (CNPq) grant 301050/2016-3 and Funda\c{c}\~ao de 
Amparo \`a Pesquisa do Estado de S\~ao Paulo (FAPESP) grant 2017/10555-0.

\end{document}